\documentclass[12pt]{amsart}
\usepackage{amsmath, amssymb, latexsym, amsthm, mathrsfs, color, mathtools, tikz,pgfplots,tikz-cd, graphicx, stackengine,url,accents, scrextend, enumitem}
\usepackage[top=2.5cm, bottom=2.5cm, left=2.5cm, right=2.5cm]{geometry}

\stackMath
\newcommand\tsup[2][2]{%
 \def\useanchorwidth{T}%
  \ifnum#1>1%
    \stackon[-1pt]{\tsup[\numexpr#1-1\relax]{#2}}{\hspace{1pt}\scriptstyle\sim}%
  \else%
    \stackon[.5pt]{#2}{\hspace{1pt}\scriptstyle\sim}%
  \fi%
}

\newcommand{\nc}{\newcommand}

\nc{\cO}{\mathcal{O}}

\newcommand{\sgg}{{\mathsf{S}_1(\Ga,\Ga)}}

\newcommand{\sgo}{{\mathsf{S}_1(\Ga,\Op)}}

\nc{\ram}{\mathsf{Ramsey}}
\nc{\soo}{\mathsf{S}_1(\cO,\cO)}
\nc{\swl}{\mathsf{S}_1(\Om,\Lambda)}
\nc{\goo}{\gone(\cO,\cO)}
\nc{\gwo}{\gone(\Om,\cO)}
\nc{\gwl}{\gone(\Om,\Lambda)}
\nc{\soox}[1]{\mathsf{S}_1(\cO(#1),\cO(#1))}
\nc{\swlx}[1]{\mathsf{S}_1(\Om(#1),\Lambda(#1))}
\nc{\goox}[1]{\gone(\cO(#1),\cO(#1))}
\nc{\gwox}[1]{\gone(\Om(#1),\cO(#1))}
\nc{\gwlx}[1]{\gone(\Om(#1),\Lambda(#1))}

\nc{\sfinoo}{\mathsf{S}_{\mathrm{fin}}(\cO,\cO)}
\nc{\sfinwl}{\mathsf{S}_{\mathrm{fin}}(\Om,\Lambda)}
\nc{\sfinww}{\mathsf{S}_{\mathrm{fin}}(\Om,\Om)}
\nc{\gfinoo}{\gfin(\cO,\cO)}
\nc{\gfinwo}{\gfin(\Om,\cO)}
\nc{\gfinwl}{\gfin(\Om,\Lambda)}
\nc{\sfinoox}[1]{\mathsf{S}_{\mathrm{fin}}(\cO(#1),\cO(#1))}
\nc{\sfinwlx}[1]{\mathsf{S}_{\mathrm{fin}}(\Om(#1),\Lambda(#1))}
\nc{\sfinwwx}[1]{\mathsf{S}_{\mathrm{fin}}(\Om(#1),\Om(#1))}
\nc{\gfinoox}[1]{\gfin(\cO(#1),\cO(#1))}
\nc{\gfinwox}[1]{\gfin(\Om(#1),\cO(#1))}
\nc{\gfinwlx}[1]{\gfin(\Om(#1),\Lambda(#1))}

\nc{\mc}{\mathcal}
\nc{\thusfar}{\my{--- Edited thus far ---}}
\nc{\lei}{\le^\oo}
\nc{\sqsubs}{\sqsubseteq^*}
\nc{\card}[1]{\left|#1\right|}
\nc{\medcard}[1]{\biggl|\,#1\,\biggr|}
\nc{\smallcard}[1]{|\,#1\,|}
\nc{\bds}{bidirectional $\roth$-scale}
\nc{\bfP}{\mathbf{P}}
\nc{\bfQ}{\mathbf{Q}}
\nc{\bbT}{\mathbb{T}}
\nc{\bbZ}{\mathbb{Z}}
\nc{\bbN}{\mathbb{N}}
\nc{\bbC}{\mathbb{C}}
\nc{\beq}{\begin{equation}}\nc{\eeq}{\end{equation}}
\nc{\mbq}{\mb{?}}
\nc{\mb}[1]{{\mbox{\textbf{#1}}}}
\nc{\nop}{$\times$}
\nc{\fbn}{\!\!\fbox{\!\nop\!}\!\!}
\nc{\yup}{\checkmark}
\nc{\forces}{\Vdash}
\nc{\name}[1]{\dot{#1}}
\nc{\tf}{\my{FINISHED THUS FAR}}
\nc{\FU}{Fr\'echet--Urysohn}
\nc{\gs}{$\gamma$~space}
\nc{\Gab}{\Gamma_{\mathrm{B}}}
\nc{\Omb}{\Omega_{\mathrm{B}}}
\nc{\Ga}{\Gamma}
\nc{\Om}{\Omega}
\nc{\smallbinom}[2]{\begin{psmallmatrix} #1\\ #2 \end{psmallmatrix}}
\nc{\bgamma}{\smallbinom{\Om}{\Ga}}

\nc{\productive}[2]{(#1,\allowbreak #2)^\x}
\nc{\prdct}[1]{#1^\x}
\nc{\Sel}{\mathsf{S}}
\nc{\sset}[2]{\{\,#1 : #2\,\}}
\nc{\smb}[1]{{\!\!\mb{#1}\!\!}}
\nc{\medset}[2]{{\biggl\{\,#1 : #2\,\biggr\}}}
\nc{\smallmedset}[2]{{\bigl\{\,#1 : #2\,\bigr\}}}
\nc{\set}[2]{{\left\{\,#1 : #2\,\right\}}}
\nc{\seq}[2]{{\la\, #1 : #2\,\ra}}
\nc{\eseq}[1]{#1_0, \allowbreak #1_1, \allowbreak\dotsc} 

\nc{\eseqint}[3]{#1_{#2}, \allowbreak\dotsc,\allowbreak #1_{#3}} 
\nc{\eseqstart}[2]{#1_{#2},\allowbreak #1_{#2+1},\dotsc } 
\nc{\eprod}[1]{#1_{1}\times \allowbreak#1_{2}\times\dotsb}

\nc{\shortprod}[1]{\prod_{n=1}^\infty{#1}_n}

\nc{\eprodint}[3]{#1_{#2}\times \allowbreak\dotsb\times\allowbreak #1_{#3}}

\nc{\seleseq}[1]{#1_1\in \mathcal{#1}_1, \allowbreak #1_2\in \mathcal{#1}_2, \allowbreak\dotsc}
\nc{\cube}{(\Cantor)^\bbN}
\nc{\Match}{\op{Match}}
\nc{\concat}[1]{\hat{\phantom{a}}\langle #1\rangle}
\nc{\poset}{\mathbb{P}}
\nc{\fn}[1]{{\op{Fn}(#1\times\w,2)}}
\nc{\linadd}{\op{linadd}}
\nc{\nonprod}{\non^\x}
\nc{\alephes}{{\aleph_0}}
\nc{\my}[1]{{\color{red}{#1}}}
\nc{\later}[1]{{\color{green} #1}}
\nc{\BTs}[1]{{\color{green} #1 (BT)}}
\nc{\Cp}{\op{C}_\mathrm{p}}
\nc{\Bp}{\op{B}_p}
\nc{\Pa}[8]{\bibitem{#1} {#2}, \emph{#3}, {#4} \textbf{#5} ({#6}), {#7}--{#8}.}
\nc{\tPa}[5]{\bibitem{#1} {#2}, \emph{#3}, {#4}, to appear.}
\nc{\sPa}[4]{\bibitem{#1} {#2}, \emph{#3}, {#4}, submitted.}
\nc{\Bc}[9]{\bibitem{#1} {#2}, \emph{#3}, in: \textbf{#4} (#5), #6 #7, #8--#9.}
\nc{\fD}{\mathfrak{D}}
\nc{\fX}{\mathfrak{X}}
\nc{\Onbd}{\Op_{\mathrm{nbd}}} 
\nc{\Omnb}{\Om_{\mathrm{nbd}}} 
\nc{\od}{\mathfrak{od}}
\nc{\Setting}[7]{\xymatrix@R=4pt@C=7pt{#1\ar@{-}[r]&#2\ar@{-}[r]&#3\\&#4\ar@{-}[u]\\
#5\ar@{-}[uu]\ar@{-}[r] & #6\ar@{-}[u]\ar@{-}[r] & #7\ar@{-}[uu]}}
\nc{\mx}[1]{\begin{matrix}#1\end{matrix}}
\nc{\plim}{p\txt{-}\lim}
\nc{\Bgp}{{\Z^\bbN}}
\nc{\Cgp}{{{\Z_2}^\bbN}}
\nc{\Cite}[1]{\textbf{[#1]}}
\nc{\Next}[1]{{#1^+}}
\nc{\cFin}{\mathrm{cF}}
\nc{\scsp}{\text{-scale space}}
\nc{\cfn}{\text{cofinal}\ }
\nc{\Con}{\text{Concentrated}}
\nc{\Lind}{\text{Lindel\"of}\,}
\nc{\con}{\text{-Concentrated}}
\nc{\lind}{\text{-Lindel\"of}\,}
\nc{\ctbl}{\text{countably }\allowbreak}
\nc{\Hur}{\text{Hurewicz}}
\nc{\intvl}[2]{{[#1(#2),\allowbreak #1(#2\!+\!1))}}
\nc{\Bdd}{\mathbf{B}}
\nc{\Dfin}{\mathfrak{D}_\mathrm{fin}}
\nc{\grbl}{{\mbox{\textit{\tiny gp}}}}
\nc{\bbP}{\mathbb{P}}
\nc{\bbH}{\mathbb{H}}
\nc{\BOfat}{\B_{\Om_{\mathrm{fat}}}}
\nc{\Bgood}{\B_{\mathrm{good}}}
\nc{\compactN}{\cl{\mathbb{N}}}
\nc{\blocks}[2]{\op{cl}_{#2}(#1)}
\nc{\blocksplus}[2]{\op{cl}^+_{#2}(#1)}
\nc{\arx}[1]{\texttt{http://arxiv.org/math/#1}}
\nc{\bq}{\begin{quote}}
\nc{\eq}{\end{quote}}
\nc{\cl}[1]{\overline{#1}}
\nc{\Cl}[2]{\mathrm{cl}_{#1}(#2)}
\nc{\CH}{the Continuum Hypothesis}
\nc{\MA}{Martin's Axiom}
\nc{\Bfat}{\B_\mathrm{fat}}
\nc{\inv}{^{-1}}
\nc{\Cantor}{{2^\w}}
\nc{\bP}{\mathbf{P}}
\nc{\bof}{\op{\fb}}
\nc{\dof}{\op{\fd}}
\nc{\bofF}{\bof(\cF)}
\nc{\sr}[3]{\underset{\mbox{#3}}{\mbox{#1}}}
\nc{\gp}{\binom{\Om}{\Ga}}
\nc{\gpsmall}{\mbox{$\gp$}}
\nc{\gig}{\gimel}
\nc{\gns}{\sone(\Om,\gig)}
\nc{\nsr}[2]{#1}
\nc{\Srg}{{\mathbb{S}}}
\nc{\Srgs}{{\mathbb{S}^*}}
\nc{\NN}{{\w^{\w}}}
\nc{\ZN}{{\Z^{\bbN}}}
\nc{\NNup}{{\bbN^{\uparrow\bbN}}}
\nc{\NNupb}{{b^{\uparrow\bbN}}}
\nc{\Pof}{\op{P}}
\nc{\PN}{{\Pof(\w)}}
\nc{\rothx}[1]{{[#1]^{\mbox{\tiny $\infty$}}}}
\nc{\tx}{{\tilde{x}}}
\nc{\roth}{{[\w]^{\w}}}
\nc{\roths}{{[b]^{\mbox{\tiny $\infty$}}}} 
\nc{\Fin}{\mathrm{Fin}}
\nc{\ici}{[\bbN]^{ \infty, \infty}}
\nc{\Inc}{{\compactN^{\uparrow\bbN}}}
\nc{\powInc}[1]{{\big(\Inc\big)^{#1}}}
\nc{\powFin}[1]{{\big(\Fin\big)^{#1}}}
\nc{\powPN}[1]{{\big(\PN\big)^{#1}}}
\nc{\NcompactN}{{\compactN^\bbN}}
\nc{\Uarrow}{\smash{\big\uparrow}}
\nc{\LE}{\preccurlyeq}
\nc{\GE}{\succcurlyeq}
\nc{\op}{\operatorname}
\nc{\im}{\op{Im}}
\nc{\Span}{\op{span}}
\nc{\maxfin}{\op{maxfin}}
\nc{\ran}{\op{range}}
\nc{\iso}{\cong}
\nc{\Madd}{{\M}^\star}
\nc{\cI}{\mathcal{I}}
\nc{\cJ}{\mathcal{J}}
\nc{\scrA}{\mathscr{A}}
\nc{\scrB}{\mathscr{B}}
\nc{\scrC}{\mathscr{C}}
\nc{\scrD}{\mathscr{D}}
\nc{\scrF}{\mathscr{F}}
\nc{\scrK}{\mathscr{K}}
\nc{\A}{\D\forall}
\nc{\B}{\mathrm{B}}
\nc{\cB}{\mathcal{B}}
\nc{\cZ}{\mathcal{Z}}
\nc{\bB}{\mathbf{B}}
\nc{\BS}{\mathbf{B}(\mathcal{S})}
\nc{\BF}{\mathbf{B}(\mathcal{F})}
\nc{\BU}{\mathbf{B}(\mathcal{U})}
\nc{\cSp}{\mathcal{S}^+}
\nc{\cFp}{\mathcal{F}^+}
\nc{\cUp}{\mathcal{U}^+}

\nc{\BG}{\B_\Ga}
\nc{\BL}{\B_\Lambda}
\nc{\BT}{\B_\Tau}
\nc{\BTstar}{\B_{\Tau^*}}
\nc{\BO}{\B_\Om}
\nc{\DO}{\cD_\Om}
\nc{\KO}{\cK_\Om}
\nc{\CG}{C_\Ga}
\nc{\CL}{C_\Lambda}
\nc{\CT}{C_\Tau}
\nc{\CTstar}{C_{\Tau^*}}
\nc{\CO}{C_\Om}
\nc{\COgp}{C_{\Om^{\grbl}}}
\nc{\CLgp}{C_{\Lambda^{\grbl}}}
\nc{\BOgp}{\B_{\Om}^{\grbl}}
\nc{\BLgp}{\B_{\Lambda^{\grbl}}}
\nc{\sfC}{\mathsf{C}}
\nc{\sfD}{\mathsf{D}}
\nc{\bD}{\mathbf{D}}
\nc{\Tau}{\mathrm{T}}
\nc{\cA}{\mathcal{A}}
\nc{\cK}{\mathcal{K}}
\nc{\cD}{\mathcal{D}}
\nc{\cF}{\mathcal{F}}
\nc{\cS}{\mathcal{S}}
\nc{\cT}{\mathcal{T}}
\nc{\cG}{\mathcal{G}}
\nc{\cY}{\mathcal{Y}}
\nc{\J}{\mathcal{J}}
\nc{\cL}{\mathcal{L}}
\nc{\cM}{\mathcal{M}}
\nc{\cN}{\mathcal{N}}
\nc{\cH}{\mathcal{H}}
\nc{\cE}{\mathcal{E}}
\nc{\One}{\mathrm{One}}
\nc{\Op}{\mathrm{O}}
\nc{\rmA}{\mathrm{A}}
\nc{\rmF}{\mathrm{F}}
\nc{\rmB}{\mathrm{B}}
\nc{\rmD}{\mathrm{D}}
\nc{\rmP}{\mathrm{P}}
\nc{\cC}{\mathcal{C}}
\nc{\cP}{\mathcal{P}}
\nc{\bbQ}{\mathbb{Q}}
\nc{\bbR}{\mathbb{R}}
\nc{\bbS}{\mathbb{S}}
\nc{\cU}{\mathcal{U}}
\nc{\cQ}{\mathcal{Q}}
\nc{\Un}{\bigcup}
\nc{\cV}{\mathcal{V}}
\nc{\cR}{\mathcal{R}}
\nc{\tcR}{\tilde{\mathcal{R}}}
\nc{\cW}{\mathcal{W}}
\nc{\Z}{{\mathbb Z}}
\nc{\Impl}{\Rightarrow}
\long\def\forget#1\forgotten{\marginpar{\textcolor{green}{Forgetting...}}}
\nc{\ft}{\mathfrak{t}}
\nc{\fb}{\mathfrak{b}}
\nc{\fc}{\mathfrak{c}}
\nc{\fd}{\mathfrak{d}}
\nc{\fg}{\mathfrak{g}}
\nc{\oo}{\infty}
\nc{\fr}{\mathfrak{r}}
\nc{\fk}{\mathfrak{k}}
\nc{\bidi}{\mathfrak{bidi}}
\nc{\fu}{\mathfrak{u}}
\nc{\fh}{\mathfrak{h}}
\nc{\fp}{\mathfrak{p}}
\nc{\fj}{\mathfrak{j}}
\nc{\fs}{\mathfrak{s}}
\nc{\w}{\omega}
\nc{\x}{\times}
\nc{\Iff}{\Leftrightarrow}

\nc{\nin}{\notin}
\nc{\cat}{\hat{\ }}
\nc{\sub}{\subseteq}
\nc{\spst}{\supseteq}
\nc{\sm}{\setminus}
\nc{\as}{\subseteq^*}
\nc{\les}{\le^*}
\nc{\leinf}{\le^{\infty}}
\nc{\leS}{\le_S}
\nc{\leF}{\le_F}
\nc{\leU}{\le_U}
\nc{\gew}{\geq^\textrm{w}}
\nc{\rest}{\restriction}

\nc{\la}{\langle}
\nc{\ra}{\rangle}
\nc{\E}{\exists}
\nc{\dom}{\op{dom}}
\nc{\cov}{\op{cov}}
\nc{\add}{\op{add}}
\nc{\addmen}{\add(\Men{})}
\nc{\cof}{\op{cof}}
\nc{\cf}{\op{cf}}
\nc{\non}{\op{non}}
\nc{\unif}{\op{non}}
\nc{\COV}{\op{COV}}
\nc{\ADD}{\op{ADD}}
\nc{\COF}{\op{COF}}
\nc{\NON}{\op{NON}}
\nc{\supp}{\op{supp}}
\nc{\impl}{\to}
\nc{\Lp}{\mathcal{L_\p}}
\nc{\Wlog}{without loss of generality}
\newtheorem{thm}{Theorem}[section]
\nc{\bthm}{\begin{thm}}
\nc{\ethm}{\end{thm}}
\newtheorem{need}[thm]{Need}
\nc{\bneed}{\begin{need}\color{dg}} \nc{\eneed}{\end{need}}
\newtheorem{prop}[thm]{Proposition}
\nc{\bprp}{\begin{prop}} \nc{\eprp}{\end{prop}}
\newtheorem{fact}[thm]{Fact}
\nc{\bfct}{\begin{fact}} \nc{\efct}{\end{fact}}
\newtheorem{prob}[thm]{Problem}
\nc{\bprb}{\begin{prob}}
\nc{\eprb}{\end{prob}}
\newtheorem{lem}[thm]{Lemma}
\nc{\blem}{\begin{lem}} \nc{\elem}{\end{lem}}
\newtheorem{app}[thm]{Application}
\nc{\bapp}{\begin{app}} \nc{\eapp}{\end{app}}
\newtheorem{claim}[thm]{Claim}
\nc{\bclm}{\begin{claim}} \nc{\eclm}{\end{claim}}
\newtheorem{cor}[thm]{Corollary}
\nc{\bcor}{\begin{cor}} \nc{\ecor}{\end{cor}}
\newtheorem{conj}[thm]{Conjecture}
\nc{\bcnj}{\begin{conj}} \nc{\ecnj}{\end{conj}}
\theoremstyle{definition}
\newtheorem{defn}[thm]{Definition}
\nc{\bdfn}{\begin{defn}} \nc{\edfn}{\end{defn}}
\newtheorem{obs}[thm]{Observation}
\nc{\bobs}{\begin{obs}} \nc{\eobs}{\end{obs}}
\theoremstyle{remark}
\newtheorem{rem}[thm]{Remark}
\nc{\brem}{\begin{rem}} \nc{\erem}{\end{rem}}
\newtheorem{cnv}[thm]{Convention}
\nc{\bcnv}{\begin{cnv}} \nc{\ecnv}{\end{cnv}}
\newtheorem{exam}[thm]{Example}
\nc{\bexm}{\begin{exam}} \nc{\eexm}{\end{exam}}
\nc{\bpf}{\begin{proof}} \nc{\epf}{\end{proof}
}
\nc{\be}{\begin{enumerate}}
\nc{\ee}{\end{enumerate}}
\nc{\bi}{\begin{itemize}}
\nc{\bimy}{\my{\begin{itemize}}
\nc{\eimy}{\end{itemize}}}
\nc{\itm}{\item}
\nc{\ei}{\end{itemize}}
\nc{\Subsection}[1]{\goodbreak\subsection*{#1}}

\nc{\sone}{\mathsf{S}_1}
\nc{\sfin}{\mathsf{S}_\mathrm{fin}}
\nc{\ufin}{\mathsf{U}_\mathrm{fin}}
\nc{\Split}{\mathsf{Split}}
\nc{\un}{\mathsf{U}_\mathrm{id}}
\nc{\unn}{\mathsf{U}_k}

\nc{\gone}{\mathsf{G}_1}    
\nc{\tgfin}{\tilde{\mathsf{G}}_\mathrm{fin}}
\nc{\gfin}{\mathsf{G}_\mathrm{fin}}

\nc{\men}[1]{\sfin(\Op(#1),\Op(#1))}
\nc{\sch}{\ufin(\cO,\Omega)}
\nc{\rothb}{\text{Rothberger}}
\nc{\pmen}{\sfin(\Omega,\Omega)}
\nc{\Rothb}{\sone(\Op,\Op)}

\nc{\prothb}{\sone(\Omega,\Omega)}
\nc{\tU}{{\tilde{U}}}
\nc{\tF}{{\tilde{F}}}
\nc{\tY}{{\tilde{Y}}}
\nc{\tX}{{\tsup[1]{X}}}
\nc{\dtX}{{\tsup[2]{X}}}
\nc{\dt}[1]{{\tsup[2]{#1}}}
\nc{\td}{{\tilde{d}}}
\nc{\tz}{{\tilde{z}}}
\nc{\cfd}{\cf(\fd)}

\nc{\msep}{\sfin(\cD,\cD)}
\nc{\rsep}{\sone(\cD,\cD)}
\nc{\cft}{\sfin(\Omega_{\mathbf{0}},\Omega_{\mathbf{0}})}
\nc{\scft}{\sone(\Omega_{\mathbf{0}},\Omega_{\mathbf{0}})}

\nc{\Umen}{U\text{-Menger}}
\nc{\hur}{\ufin(\cO,\Gamma)}
\nc{\tUmen}{\tU\text{-Menger}}

\nc{\Men}{\text{Menger}}
\nc{\Sch}{\text{Scheepers}}

\nc{\aspst}{\prescript{*}{}{\spst}\ }
\nc{\eqs}{=^*}
\nc{\ctblOm}{\Omega_{\mathrm{ctbl}}}
\nc{\GNga}{{\smallbinom{\Om}{\Ga}}}
\nc{\ctblga}{\smallbinom{\ctblOm}{\Ga}}

\nc{\nadd}{\cN_{\mathrm{add}}}
\nc{\ball}{\mathrm{B}}
\nc{\cOunif}{\cO^{\textrm{unif}}}

\nc{\FS}{\op{FS}}
\nc{\sums}{\op{SS}}
\nc{\SG}{\op{SG}}
\nc{\tSG}{\op{\widetilde{SG}}}
\nc{\G}{\op{G}}
\nc{\pBM}{\op{pBM}}
\nc{\FSG}{\op{FSG}}
\nc{\M}{\op{M}}
\nc{\FP}{\op{FP}}
\nc{\nonNadd}{\non(\nadd)}
\nc{\borga}{\Ga_\mathrm{Bor}}
\nc{\pick}{x}
\nc{\gen}{y}
\nc{\nullzind}{\sone(\{\Op_n^{\mathsf{unif}}\}_{n\in\bbN},\Ga)}
\nc{\nullzindf}[1]{\sone(\{\Op_{#1}^{\mathsf{unif}}\}_{n\in\bbN},\Ga)}

\definecolor{dg}{RGB}{42,101,24}

\nc{\myb}[1]{\textcolor{blue}{#1}}
\nc{\mydg}[1]{{\color{dg}{#1}}}
\DeclareMathOperator{\eexists}{\exists}
\DeclareMathOperator{\fforall}{\forall}
\nc{\Exists}[1]{\bigl(\eexists #1\bigr)}
\nc{\Forall}[1]{\bigl(\fforall #1\bigr)}
\nc{\End}[1]{\bigl(#1\bigr)}
\nc{\dmo}[2]{\DeclareMathOperator{#1}{#2}}
\dmo{\Asc}{Asc}
\nc{\plusmin}{\wedge}

\nc{\cBsub}{{\cB^{\mbox{\tiny $\sub$}}}}

\nc{\Alice}{{\textsc{Alice}{}}}
\nc{\Bob}{{\textsc{Bob}}}
\nc{\BM}{\op{BM}}
\nc{\Palpha}{{\bbS_\alpha}}
\nc{\Pbeta}{{\bbS_\beta}}
\nc{\Pwtwo}{\bbS_{\w_2}}
\nc{\restrict}{\upharpoonright}

\usepackage{float}
\usepackage[all]{xy}

\author[P. Szewczak]{Piotr Szewczak}
\address{Piotr Szewczak, 
 Institute of Mathematics, Faculty of Mathematics, Informatics, and Mechanics, University of Warsaw, Banacha 2, 02-097, Warsaw, Poland
}
\email{p.szewczak@wp.pl}
\urladdr{https://piotrszewczak.pl}

\author[T.~Weiss]{Tomasz Weiss}\address{Tomasz Weiss, Institute of Mathematics, Faculty of Mathematics and Natural Science College of Sciences, Cardinal Stefan Wyszy\'nski University in Warsaw, W\'oycickiego 1$\slash$3, 01--938 Warsaw, Poland}
\email{tomaszweiss@o2.pl}

\author[L.~Zdomskyy]{Lyubomyr Zdomskyy}
\address{Lyubumoyr Zdomskyy, Institut f\"ur Diskrete Mathematik und Geometrie, Technische Universit\"at Wien, Wiedner Hauptstrasse 8-10/104, 1040 Wien, Austria.}
\email{lzdomsky@gmail.com}
\urladdr{https://dmg.tuwien.ac.at/zdomskyy/}

\title{Small Hurewicz and Menger sets which have large continuous images}

\subjclass[2010]{Primary: 54D20 
}
\keywords{Hurewicz, Menger, totally imperfect, perfectly meager, everywhere meager, countably perfectly meager}
\thanks{
The research of the first and second authors was funded by the National Science Center, Poland Weave-UNISONO call in the Weave programme
Project: Set-theoretic aspects of topological selections 2021/03/Y/ST1/00122. 
The research of the third author
was funded in whole by the Austrian Science Fund
(FWF) [10.55776/I5930 and 10.55776/PAT5730424].
}

\begin{document}

\maketitle
\begin{abstract}
We provide new techniques to construct sets of reals without perfect subsets and with the Hurewicz or Menger covering properties.
In particular, we show that if \CH{} holds, then there are such sets which can be mapped continuously onto the Cantor space.
These results allow to separate the properties of Menger and $\sgo$ in the realm of sets of reals without perfect subsets and solve a problem of Nowik and Tsaban concerning perfectly meager subsets in the transitive sense.
We present also some other applications of the mentioned above methods.
\end{abstract}

\section{Introduction}

Throughout this paper, we assume that a set $X$ with a given topological property is a subset of the Cantor space $\Cantor$ or its square.
A set $X$ is \emph{Menger} if for any sequence $\eseq{\cU}$ of open covers of $X$, there are finite sets $\cF_0\sub\cU_0, \cF_1\sub\cU_1,\dotsc$ such that the family $\sset{\Un\cF_n}{n\in\w}$ covers $X$.
If in addition the family $\sset{\Un\cF_n}{n\in\w}$ is a \emph{$\gamma$-cover} of $X$, i.e.,  for every $x\in X$, the set $\sset{n}{x\in\Un\cF_n}$ is co-finite, then the set $X$ is \emph{Hurewicz}.
We have the following implications between considered properties.
\[
\sigma\text{-compact}\longrightarrow \text{Hurewicz}\longrightarrow \text{Menger}.
\]
A set is \emph{totally imperfect} if it does not contain a perfect subset, i.e., a homeomorphic copy of $\Cantor$.
By the results of Bartoszy\'{n}ski--Shelah~\cite{BaSh01} and Tsaban--Zdomskyy~\cite{sfh}, there are uniform constructions of totally imperfect sets in ZFC which show that none of the above implications are reversible in the realm of sets of reals.
The above properties are central in the combinatorial covering properties theory and found applications in such areas as forcing~\cite{ChoZdo}, function spaces~\cite{coc6}, Ramsey theory in algebra~\cite{Tsa18}, combinatorics of discrete subspaces~\cite{AurD}, hyperspaces with the Vietoris topology~\cite{krup}, products of Lindel\"of spaces~\cite{AurTall, pMGen} and products of paracompact spaces~\cite{Sze17}.

Each of the considered properties has a combinatorial characterization using continuous mappings into the Baire space $\NN$.
A set $X$ is Menger (Hurewicz) if and only if every continuous image of $X$ into $\NN$ is non-dominating (bounded)~\cite{reclaw}.
It follows that each of these properties has a \emph{critical cardinality}, i.e., the minimal cardinality of a subset of $\NN$ not having a given property.
The critical cardinality for Menger is the dominating number $\fd$ and for Hurewicz it is the bounding number $\fb$.

Many known constructions of totally imperfect sets which are Menger or Hurewicz involve category-theoretic, measure-theoretic, combinatorial or forcing methods.
However, each of these constructions provide a set satisfying $\sgo$, which is another classical covering property well studied in the field~\cite{coc2, sakai, SPMProd}.
A set $X$ is $\sgo$ if for every sequence $\eseq{\cU}$ of open $\gamma$-covers of $X$, there are sets $U_0\in\cU_0, U_1\in\cU_1,\dotsc$ such that the family $\sset{U_n}{n\in\w}$ covers $X$.
The critical cardinality for $\sgo$ is $\fd$ and the following implications hold~\cite{coc2}. 

\begin{figure}[H]
\begin{center}
\begin{tikzcd}[ampersand replacement=\&,column sep=1cm]
{} \&  \text{perfectly meager}  \& \& \\
\sigma\text{-compact}\arrow{r}\& \text{Hurewicz}\arrow{r}\arrow[u, left,"+ \text{ totally imperfect}" swap]\&\text{Menger}\&{}\\
{}\& {}\& \sgo\arrow{u}\arrow{r} \& \text{ totally imperfect}
\end{tikzcd}
\end{center}
\end{figure}
The only ZFC examples for separation Hurewicz or Menger from $\sgo$ are not satisfactory, i.e., they are the Cantor space and, more generally, sets containing perfect subsets.
After submitting the paper, Zakrzewski shared with us earlier hand notes of Pol with similar constructions as presented in our Theorems~\ref{thm:main} and \ref{thm:mainmen}, where in particular we construct a totally imperfect Hurewicz set and totally imperfect Menger set, using additional set-theoretic assumptions beyond ZFC.   
Professor Roman Pol also pointed out to us that our Theorem~\ref{thm:mainmen} can be deduced from papers~\cite{Pol,PPol} written in a different context and using different techniques.

A set $X\sub\Cantor$ is \emph{perfectly meager} if for any perfect set $P\sub \Cantor$, the set $X\cap P$ is meager in the relative topology of $P$.
By the result of Just--Miller--Scheepers--Szeptycki~\cite[Theorem~5.5]{coc2} any totally imperfect Hurewicz subset of $\Cantor$ is perfectly meager.
It follows from the result of Bartoszy\'{n}ski--Shelah~\cite{BaSh01} that in ZFC there is a perfectly meager set which is not Hurewicz.
Let $\cof(\cM)$ be the minimal cardinality of a family of meager sets in $\Cantor$ such that each meager set in $\Cantor$ is contained in a memebr of the family.
Using techinques from the work of Szewczak--Wiśniewski~\cite[Proposition~2.9.]{Lusin}, it can be shown that assuming the equality $\fd=\cof(\cM)$, there is a totally imperfect Menger set which is not perfectly meager.

Using various set-theoretic assumptions which are weak portions of \CH{}, we construct totally imperfect Hurewicz or Menger sets which can be mapped continuously onto the Cantor cube.
Most of already known strategies to construct totally imperfect sets with these properties have substantial limitations which prevent from reaching such a goal.
These results have consequences for separation of considered properties.
Another benefit is a solution of the below problem which was formulated independently by Nowik and Tsaban (in oral communication).
By $\oplus$ we mean a standard group operation in $\Cantor$, i.e., for every $x,y\in\Cantor$ we have 
\[
(x\oplus y)(n):=(x(n)+y(n)) (\bmod 2)
\]
for all natural numbers $n$.
For a set $A\sub\Cantor$ and an element $t\in\Cantor$, let $t\oplus A:=\sset{t\oplus a}{a\in A}$.
A set $X\sub\Cantor$ is \emph{perfectly meager in the transitive sense}~\cite{nsw, NW} (PMT in short) if for any perfect set $P\sub\Cantor$, there is a $\sigma$-compact set $F\sub\Cantor$ containing $X$ such that the set $F\cap (t\oplus P)$ is meager in the relative topology of $t\oplus P$ for all $t\in\Cantor$. 

\bthm[{Nowik~\cite[Theorem~1.]{nowik}}]\label{thm:nowik}
Each Hurewicz set in $\Cantor$ which cannot be mapped continuously onto $\Cantor$ is PMT.
\ethm

In view of Zakrzewski's result~\cite{zakp} which states that every totally imperfect Hurewicz set is universally meager (see definition in Subsection~\ref{sec:EM}), the following question seems to be natural.

\bprb[Nowik, Tsaban]\label{prb:pmt}
Is any Hurewicz totally imperfect subset of $\Cantor$, PMT? 
\eprb

We show that a solution to this problem is independent from ZFC.
\section{Hurewicz sets}

Let $\cov(\cN)$ be the minimal cardinality of a family of Lebesgue-null subsets of $\Cantor$ which covers $\Cantor$.
A subset $X$ of a topological space is a \emph{$\lambda'$-set} if for any countable subset $A$ of the space the set $A$ is a $G_\delta$-set in the relative topology of $X\cup A$.

\bthm\label{thm:main}
Assume that $\cov(\cN)=\fb=\fc$.
Then there exists a set $X\sub\Cantor\x\Cantor$ with the following properties:
\be
\item $X$ is Hurewicz,
\item $X$ is a $\lambda'$-set,
\item $X$ is totally imperfect,
\item $\pi[X]=\Cantor$,
\item a homeomorphic copy of $X$ in $\Cantor$ is not PMT.
\ee
\ethm

\brem
By the equivalence (1) $\Leftrightarrow$ (3) in Proposition 2.6 from~\cite{a}, in the Laver model every $\lambda'$-set is Hurewicz, and by the discussion before Lemma~2.3 in \cite{b} we have that every Hurewicz $\lambda$-set in this model has size $\w_1$.
Since $\fb=\fc=\w_2$ in the Laver model, we conclude that $\fb=\fc$ is not sufficient to achieve the disjunction of (2) and (4) in Theorem~\ref{thm:main}.
\erem

In order to prove Theorem~\ref{thm:main} we need the following auxiliary results.

\blem[{\cite[Lemma~6]{NW}}]\label{lem:NW}
A set $X\sub \Cantor$ is PMT if and only if for any sequence $\eseq{P}$ of perfect sets in $\Cantor$, there are compact sets $\eseq{F}\sub \Cantor$ such that $X\sub \Un_{n\in\w}F_n$ and $t\oplus P_n\not\sub F_m$ for every natural numbers $n,m$ and element $t\in\Cantor$.
\elem

For every set $A\sub \Cantor\x\Cantor$ and element $x\in 2^\w$, let 
\[
A\vert_x:=\sset{y\in \Cantor}{(x,y)\in A}.
\]
Let $\pi\colon\Cantor\x\Cantor\to\Cantor$ be the projection onto the first coordinate.

\blem\label{lem:step1}
Let $A\sub\Cantor\x\Cantor$ be a clopen set such that $\pi[A]=\Cantor$ and $M\sub\Cantor\x\Cantor$ be a closed set such that the sets $M\vert_x$ are nowhere dense in $\Cantor$ for all elements $x\in\Cantor$.
Then there is a clopen set $B\sub A\sm M$ such that $\pi[B]=\Cantor$.
\elem

\bpf
Any set $A$ with the given property can be replaced by a finite union of nonempty basic open sets in $\Cantor\x\Cantor$.
Since the sets $M\vert_x$ are nowhere dense in $\Cantor$ for all $x\in\Cantor$, it is enough to prove the statement for $A=\Cantor\x\Cantor$.
For each element $x\in \Cantor$ there are an element $y\in\Cantor$ and clopen basic sets $U_x,V_x\sub \Cantor$ such that 
\[
(x,y)\in U_x\x V_x\sub (\Cantor\x\Cantor) \sm M.
\]
Then there is a finite set $X'\sub X$ such that the family $\sset{U_x}{x\in X'}$ covers $\Cantor$.
Put
\[
B:=\Un\sset{U_x\x V_x}{x\in X'}.\qedhere
\]
\epf

\blem\label{lem:step2}
Let $\eseq{M}$ be an increasing sequence of closed sets in $\Cantor\x\Cantor$ such that the sets $M_n\vert_x$ are nowhere dense in $\Cantor$ for all natural numbers $n$ and every $x\in \Cantor$.
Then there is a compact set $D\sub\Cantor\x\Cantor$ disjoint from $\Un_{n\in\w}M_n$ such that $\pi[D]=\Cantor$.
\elem

\bpf
Let $D_0:=\Cantor\x\Cantor$.
Fix a natural number $n$ and assume that a clopen set $D_n\sub\Cantor\x\Cantor$ with $D_n\cap M_n =\emptyset$ and $\pi[D_n]=\Cantor$, has already been defined.
Let $D_{n+1}$ be a set from Lemma~\ref{lem:step1}, applied to the sets $D_n$ and $M_n$.
Then the set $D:=\bigcap_{n\in\w}D_n$ is disjoint from $\Un_{n\in\w}M_n$ and $\pi[D]=\Cantor$.
\epf

\blem\label{lem:nonPMT}
Let $X\sub\Cantor\x\Cantor$ be a set such that for each open set $G\sub\Cantor\x\Cantor$ containing $X$ there is an element $x\in\Cantor$ with $G\vert_x=\Cantor$.
Then there is a homeomorphic copy of $X$ in $\Cantor$ which is not PMT.
\elem

\bpf
Let $K_0,K_1\sub \Cantor$ be perfect sets such that the map $\oplus\colon K_0\x K_1\to \Cantor$ is injective and let $\phi_0,\phi_1$ be homeomorphisms from $\Cantor$ onto $K_0$ and $K_1$, respectively.
Then the map 
\[
\Phi\colon \Cantor\x\Cantor\to K_0\oplus K_1
\]
 defined by 
\[
\Phi(x,y):=\phi_0(x)\oplus \phi_1(y),
\]
for all $(x,y)\in \Cantor\x\Cantor$, is a homeomorphism.

Assume that the set $\Phi[X]$ is PMT.
Let $\sset{P_n}{n\in\w}$ be an enumeration of a clopen basis of $K_1$.
Applying Lemma~\ref{lem:NW} to the sequence $\eseq{P}$ and to the set $\Phi[X]$, there are compact sets $\eseq{F}\sub\Cantor$ such that $\Phi[X]\sub\Un_{n\in\w}F_n$ and $t\oplus P_n\not\sub F_m$ for every natural numbers $n,m$ and elements $t\in\Cantor$.
Consequently, the sets $F_n$ are nowhere dense in $t\oplus K_1$ for all $t\in \Cantor$.
Let $M_n:=\Phi\inv[\Un_{i\leq n}F_i]$ for all natural numbers $n$.
Then $X\sub \Un_{n\in\w}M_n$ and the sets $M_n\vert_x$ are nowhere dense in $\Cantor$ for every natural number $n$ and each $x\in\Cantor$.
Let $D\sub\Cantor\x\Cantor$ be a set obtained from Lemma~\ref{lem:step2}.
Then the open set $O:=(\Cantor\x\Cantor)\sm A$ contains $X$ and $O\vert_x\neq\Cantor$ for all elements $x\in\Cantor$, which is a contradiction. Thus, the set $\Phi[X]$ is not PMT.
\epf

Let $\mu$ be the Lebesgue measure on $\Cantor$.

\blem \label{l1}
Let $O\sub \Cantor\x\Cantor$ be an open set 
such that $\mu(O\vert_x)=1$ for all elements $x\in \Cantor$.
Then for every real number $a\in (0,1)$ there exists 
a compact set $K\sub O$ such that $\mu(K\vert_x)>1-a$ for all elements $x\in \Cantor$.
\elem

\bpf
For each $x\in \Cantor$ there is an an open set $U_x\sub\Cantor$ such that $\Cantor\sm O\vert_x\sub U_x$ and $\mu(U_x)<a$.
Then there is a clopen neighborhood $V_x\sub\Cantor$ of $x$ such that
$\Cantor\sm O\vert_{x'}\sub U_x$ for all $x'\in V_x$.
Let $\{V_{x_i}:i< n\}$ be a finite subcover of the cover $\{V_x:x\in 2^\w\}$
of $2^\w$.
Replacing $V_{x_i}$ with $V_{x_i}\setminus\bigcup_{j<i}V_{x_j}$,
if necessary, we may assume in addition that the sets $V_{x_i}$
are disjoint.
Define 
\[
W:=\bigcup_{i<n}V_{x_i}\times U_{x_i}.
\]
Then 
\[
(\Cantor\x\Cantor)\sm O\sub W\text{\quad and\quad
}\mu(W\vert_x)<a
\]
for all $x\in 2^\w$.
Take 
\[
K:=(\Cantor\x\Cantor)\sm W.\qedhere
\]
\epf

\bcor \label{c1}
Let $G\sub \Cantor\x\Cantor$ be a $G_\delta$-set 
such that $\mu(G\vert_x)=1$ for all elements $x\in 2^\w$.
Then for every real number $a\in (0,1)$, there exists 
a compact set $K\sub G$ such that $\mu(K\vert_x)>1-a$ for all elements $x\in 2^\w$.
Consequently, there exists a $\sigma$-compact set $F\sub G$
such that $\mu(F\vert_x)=1$ for all $x\in \Cantor$.
\ecor

\bpf
Let us write $G$ in the form $\bigcap_{n\in\w} O_n$ where each set $O_n$ is open and $O_{n+1}\sub O_n$ for all natural numbers $n$.
By Lemma~\ref{l1} for each natural number $n$ there exists a compact set $K_n\sub O_n$ such that $\mu(K_n\vert_x)>1-\frac{a}{2^{n+1}}$ for all $x\in 2^\w$.
Take $K:=\bigcap_{n\in\w}K_n$.
\epf

\bcor \label{cor:c2}
Let $G\sub \Cantor\x\Cantor$ be a $G_\delta$-set 
such that $\mu(G\vert_x)=1$ for all elements $x\in 2^\w$ and $\mu(\sset{x\in\Cantor}{G\vert_x=\Cantor})=1$.
Then there exists a $\sigma$-compact set $F\sub G$
such that $\mu(F\vert_x)=1$ for all elements $x\in \Cantor$ and $\mu(\sset{x\in\Cantor}{F\vert_x=\Cantor})=1$.
\ecor

\bpf
Let $F'\sub G$ be a $\sigma$-compact set from Corollary~\ref{c1}.
There is an $F_\sigma$-set $A\sub \sset{x\in\Cantor}{G\vert_x=\Cantor}$ with $\mu(A)=1$. 
Take $F:=F'\cup (A\x\Cantor)$.
\epf

The following topological characterization of Hurewicz sets due to Just--Miller--Scheepers--Szeptycki plays a crucial role in our considerations.

\blem[{\cite[Theorem~5.7.]{coc2}}]\label{lem:HurFsigma}
A subset $X$ of a compact space $K$ is Hurewicz if and only if for any $G_\delta$-set $G\sub K$ containing $X$, there is an $F_\sigma$-set set $F\sub K$ such that $X\sub F\sub G$.
\elem

\bpf[{Proof of Theorem~\ref{thm:main}}]
Let $\cG$ be the family of all $G_\delta$-sets $G$ in $\Cantor\x\Cantor$ such that $\mu(G\vert_x)=1$ for all $x\in\Cantor$ and $\mu(\sset{x\in\Cantor}{G\vert_x=\Cantor})=1$, 
$\cH$ be the family of all other $G_\delta$-sets in $\Cantor\x\Cantor$ and let $\cF$ be the family of all $F_\sigma$-sets
$F$ in $\Cantor\x\Cantor$ such that $\mu(F\vert_x)=1$ for all $x\in \Cantor$ and $\mu(\sset{x\in\Cantor}{F\vert_x=\Cantor})=1$.
Suppose that $\sset{(G_\alpha,H_\alpha)}{\alpha<\fc}$ is an enumeration of sets in $\cG\x\cH$ and $2^\w=\{x_\alpha:\alpha<\fc\}$.
By recursion on $\alpha<\fc$, we construct sets $F_\alpha\in\cF$ and elements $u_\alpha, w_\alpha\in \Cantor\x\Cantor$ with the following properties:

\be[label=(\it\roman*)]
\item $\sset{u_\beta, w_\beta}{\beta<\alpha}\cap G_\alpha\sub F_\alpha\sub G_\alpha$, \label{itm:1hur}
\item $u_\alpha,w_\alpha\in \bigcap\sset{F_\beta}{\beta\leq \alpha}$, \label{itm:2hur}
\item \label{itm:Halpha} $w_\alpha\notin H_\alpha$,
\item $\pi(u_\alpha)=x_\alpha$.\label{itm:4hur}
\ee 

By Corollary~\ref{cor:c2} there is a  set $F_0\in\cF$ with $F_0\sub G_0$.
Since $H_0\notin\cG$, we have
\[
\mu(\sset{x\in\Cantor}{F_0\vert_x=\Cantor})=1
\text{\quad and \quad}
\mu(\sset{x\in\Cantor}{H_0\vert_x\neq\Cantor})>0,
\]
 or there is an element $x\in\Cantor$ such that 
 \[
\mu(F_0\vert_x)=1\text{\quad and \quad}\mu(\Cantor\sm H_0\vert_x)>0.
 \]
In both cases, there is an element $x\in\Cantor$ such that 
\[
\bigl(F_0\vert_x\bigr)\sm\bigl( H_0\vert_x\bigr)\neq\emptyset.
\]
Pick elements
\[
u_0\in F_0\cap\bigl(\{x_0\}\x\Cantor\bigr)\text{\quad and\quad}
w_0\in F_0\cap \bigl(\{x\}\x\Cantor\bigr)\sm H_0.
\]

Fix an ordinal number $\alpha<\fc$ and assume that the objects mentioned above have been constructed for all $\beta<\alpha$.
The set $\sset{u_\beta, w_\beta}{\beta<\alpha}\cap G_\alpha$ has cardinality smaller than $\fb$, which is the critical cardinal number for the Hurewicz property, and thus this set is Hurewicz.
By Lemma~\ref{lem:HurFsigma} and Corollary~\ref{cor:c2} there is a set $F_\alpha\in\cF$ such that 
\[
\sset{u_\beta,w_\beta}{\beta<\alpha}\cap G_\alpha\sub F_\alpha\sub G_\alpha.
\] 
The intersection of less than $\cov(\cN)$ subsets of $\Cantor$, where each subset has measure $1$, and a Borel set with positive measure, is nonempty.
Since $H_\alpha\notin\cG$ for all  $\beta\leq \alpha$ we have
\[
\mu(\sset{x\in\Cantor}{F_\beta\vert_x=\Cantor})=1
\text{\quad and \quad}
\mu(\sset{x\in\Cantor}{H_\alpha\vert_x\neq\Cantor})>0
\]
 or there is an element $x\in\Cantor$ such that for all $\beta\leq\alpha$ we have
 \[
\mu(F_\beta\vert_x)=1\text{\quad and \quad}\mu(\Cantor\sm H_\alpha\vert_x)>0.
 \]
Since $\alpha<\cov(\cN)$, in both cases there is an element $x\in\Cantor$ such that  
\[
\bigl(\bigcap_{\beta\leq \alpha}F_\beta\vert_{x}\bigr)\sm\bigl( H_\alpha\vert_x\bigr)\neq\emptyset.
\]
Pick elements
\[
u_\alpha\in \bigcap_{\beta\leq\alpha}F_\beta\cap\bigl(\{x_\alpha\}\x\Cantor\bigr)
\text{\quad and\quad }
w_\alpha\in \Bigl( \bigcap_{\beta\leq\alpha}F_\beta\cap \bigl(\{x\}\x \Cantor\bigr)\Bigr)\sm H_\alpha.
\]
This completes our recursive construction.
Let $X:=\sset{u_\alpha,w_\alpha}{\alpha<\fc}$.

(1)
Use a characterization from Lemma~\ref{lem:HurFsigma}.
Let $G\sub\Cantor\x\Cantor$ be a $G_\delta$-set containing $X$.
By~\ref{itm:Halpha}, none of the sets $H_\alpha$ contains $X$.
Thus, there is an ordinal number $\alpha<\fc$ such that $G=G_\alpha$.
By~\ref{itm:1hur} and~\ref{itm:2hur}, we have $X\sub F_\alpha\sub G_\alpha$.

(2)
Let $A$ be a countable subset of $\Cantor\x\Cantor$.
Then there is an ordinal number $\xi$ such that $G_\xi=(\Cantor\x\Cantor)\sm A$.
By~(1) and~(2), we have 
\[
X\sm A=X\cap G_\xi\sub F_\xi\sub G_\xi.
\]
Thus, the set $A$ is a $G_\delta$-set in the relative topology of $X\cup A$.

(3) Assume that $X$ contains a perfect set $P$.
Then there is a countable set $A\sub P$ such that the set $P\sm A$ is homeomorphic with the Baire space $\NN$.
Since the set $X$ is a $\lambda'$-set, the set $P\sm A$ is an $F_\sigma$-subset of $X$.
Since the Hurewicz property is preserved by $F_\sigma$-subsets, the set $P\sm A$ is Hurewicz too, which is a contradiction.

(4) It follows from~\ref{itm:4hur}.

(5) Apply Lemma~\ref{lem:nonPMT}.
\epf

\bcor
Assume that $\cov(\cN)=\fb=\fc$.
There is a Hurewicz totally imperfect set which is not $\sgo$.
\ecor

\brem
The proof of Theorem~\ref{thm:main} actually shows that for the minimal cardinality $\kappa$ of a subset of $\Cantor$ which cannot be covered by an $F_\sigma$ measure zero set and any set $Y\sub\Cantor\x\Cantor$ with $\card{Y}<\kappa$, the set $Y$ is a $G_\delta$-subset of $Y\cup X$.
\erem
\section{Menger sets}

Let $\cov(\cM)$ be the minimal cardinality of a family of meager subsets of $\Cantor$ which covers $\Cantor$.

\bthm\label{thm:mainmen}
Assume that $\cov(\cM)=\fc$.
Then there is a Menger totally imperfect set $X\sub\Cantor\x\Cantor$ with $\pi[X]=\Cantor$.
\ethm

\blem\label{lem:step3}
Let $\eseq{M}$ be an increasing sequence of closed sets in $\Cantor\x\Cantor$ such that the sets $M_n\vert_x$ are nowhere dense in $\Cantor$ for all natural numbers $n$ and every $x\in \Cantor$.
Then there is a $\sigma$-compact set $D\sub\Cantor\x\Cantor$ disjoint from $\Un_{n\in\w}M_n$ such that the sets $D\vert_x$ are dense in $\Cantor$ for all elements $x\in\Cantor$.
\elem

\bpf
Let $\eseq{B}$ be an enumeration of all basic clopen sets in $\Cantor\x\Cantor$.
Fix a natural number $m$.
By Lemma~\ref{lem:step2}, there is a compact set $D_m\sub B_m$ disjoint from the set $\Un_{n\in\w}M_n$ such that $\pi[D_m]=\pi[B_m]$.
Take $D:=\Un_{m\in\w}D_m$.
\epf

Let $\cG$ be the family of all $G_\delta$-sets $G\sub\Cantor\x\Cantor$ such that the sets $G\vert_x$ are dense in $\Cantor$ for all elements $x\in\Cantor$.

\blem\label{lem:men}
Let $X\sub\Cantor\x\Cantor$ be a Menger set and $\eseq{\cU}$ be a sequence of open covers of $X$, where each family $\cU_n$ is open in $\Cantor\x\Cantor$ and the set $\bigcap_{n\in\w}\Un\cU_n$ is in $\cG$.
Then there are finite sets $\cF_0\sub\cU_0, \cF_1\sub\cU_1,\dotsc$ such that the family $\Un_{n\in\w}\cF_n$ covers $X$ and the set $\Un_{n\in\w}\Un\cF_n$ is in $\cG$.  
\elem

\bpf
Let $M_n:=\Cantor\x\Cantor\sm\Un\cU_n$ for all natural numbers $n$.
Let $D\sub\Cantor$ be a $\sigma$-compact set from Lemma~\ref{lem:step3} applied to the sequence $\eseq{M}$.
Since the set $X\cup D$ is Menger, there are finite sets $\cF_n\sub\cU_n$ for all natural numbers $n$ such that the family $\Un_{n\in\w}\cF_n$ covers $X\cup D$.
\epf

\bpf[{Proof of Theorem~\ref{thm:mainmen}}]
Let $\smallmedset{(\cU_0^{(\alpha)}, \cU_1^{(\alpha)}, \dotsc)}{\alpha<\fc}$ be an enumeration of all countable sequences $\eseq{\cU}$ of open families in $\Cantor\x\Cantor$ such that the sets $\bigcap_{n\in\w}\Un\cU_n$ are in $\cG$.
Let $\smallmedset{(\cV_0^{(\alpha)}, \cV_1^{(\alpha)}, \dotsc)}{\alpha<\fc}$ be an enumeration of all other countable sequences of open families in $\Cantor\x\Cantor$, let $\sset{P_\alpha}{\alpha<\fc}$ be an enumeration of all perfect sets in $\Cantor\x\Cantor$ and $\Cantor=\sset{x_\alpha}{\alpha<\fc}$.

By recursion on $\alpha<\fc$, we construct finite sets  $\cF_0^{(\alpha)}\sub \cU_0^{(\alpha)}, \cF_1^{(\alpha)}\sub \cU_1^{(\alpha)}, \dotsc$, elements $u_\alpha,w_\alpha\in \Cantor\x\Cantor$ and $p_\alpha\in P_\alpha$ with the following properties:

\be[label=(\it\roman*)]
\item the set $F_\alpha:=\Un_{n\in\w}\Un\cF_n^{(\alpha)}$ is in $\cG$,
\item $\sset{u_\beta,w_\beta}{\beta<\alpha}\cap G_\alpha\sub F_\alpha$, where $G_\alpha:=\bigcap_{n\in\w}\Un\cU^{(\alpha)}_n$ ,\label{itm:2men}
\item $u_\alpha,w_\alpha \in\bigcap_{\beta\leq\alpha}F_\beta$,\label{itm:3men}
\item $u_\alpha, w_\alpha\notin \sset{p_\beta}{\beta\leq \alpha}$,\label{itm:4men}
\item $w_\alpha\notin H_\alpha$, where $H_\alpha:=\bigcap_{n\in\w}\Un\cV_n^{(\alpha)}$,\label{itm:5men}
\item $\pi(w_\alpha)=x_\alpha$.
\ee

Pick $p_0\in P_0$.
By Lemma~\ref{lem:men} there are finite sets $\cF_0^{(0)}\sub \cU_0^{(0)}, \cF_1^{(0)}\sub \cU_1^{(0)}, \dotsc$ such that the set $F_0:=\Un_{n\in\w}\Un\cF_n^{(0)}$ is in $\cG$.
The sets $F_0\vert_{x}$ are dense in $\Cantor$ for all elements $x\in\Cantor$.
Since the set $H_0:=\bigcap_{n\in\w}\Un\cV_n^{(0)}$ is not in $\cG$, there is an element $x\in\Cantor$ such that the set $H_0\vert_x$ is not dense in $\Cantor$, and thus
\[
\bigl(F_0\vert_x\bigr)\sm\bigl( H_0\vert_x\bigr)\neq\emptyset.
\]
Pick elements 
\[
u_0\in\bigl(F_0\cap (\{x_0\}\x\Cantor)\bigr)\sm\{p_0\}\text{\quad and\quad }
w_0\in F_0\sm\bigl(H_0\cup\{p_0\}\bigr).
\]

Fix an ordinal number $\alpha<\fc$ and assume that the objects mentioned above have been constructed for all $\beta<\alpha$.
Pick $p_\alpha\in P_\alpha\sm\sset{u_\beta,w_\beta}{\beta<\alpha}$.
Since $\alpha<\cov(\cM)\leq\fd$, the set $\sset{u_\beta,w_\beta}{\beta<\alpha}\cap G_\alpha$ has cardinality smaller than $\fd$, which is the critical cardinal number for the Menger property, and thus this set is Menger. 
By Lemma~\ref{lem:men} there are finite sets $\cF_0^{(\alpha)}\sub \cU_0^{(\alpha)}, \cF_1^{(\alpha)}\sub \cU_1^{(\alpha)}, \dotsc$ such that the set $F_\alpha:=\Un_{n\in\w}\Un\cF_n^{(\alpha)}$ is in $\cG$ and 
\[
\sset{u_\beta,w_\beta}{\beta<\alpha}\cap G_\alpha\sub F_\alpha.
\]
The intersection of less than $\cov(\cM)$
dense $G_\delta$-sets in $\Cantor$ and a subset of $\Cantor$ with nonempty interior is nonempty.
For each $\beta\leq\alpha$
the sets $F_\beta\vert_x$ are dense $G_\delta$-sets in $\Cantor$ for all elements $x\in\Cantor$.
Since the set $H_\alpha:=\bigcap_{n\in\w}\Un\cV^{(\alpha)}_n$ is not in $\cG$, there is $x\in\Cantor$ such that the set $H_\alpha\vert_x$ is not dense in $\Cantor$, i.e., the set $\Cantor\sm H_\alpha\vert_x$ has nonempty interior.
Since $\alpha<\cov(\cM)$, we can pick elements
\[
u_\alpha\in\bigl( \bigcap_{\beta\leq\alpha}F_\beta\cap(\{x_\alpha\}\x\Cantor)\bigr)\sm \sset{p_\beta}{\beta\leq \alpha}
\]
and
\[
w_\alpha\in \bigcap_{\beta\leq\alpha}F_\beta\sm\bigl(H_0\cup\sset{p_\beta}{\beta\leq \alpha}\bigr)
\]
This completes our recursive construction.
Let $X:=\sset{u_\alpha,w_\alpha}{\alpha<\fc}$.

Let $\eseq{\cU}$ be a sequence of open covers of $X$.
For each ordinal number $\alpha<\fc$ we have $w_\alpha\notin H_\alpha=\bigcap_{n\in\w}\Un\cV_n^{(\alpha)}$, and thus none of the sequences $\cV_0^{(\alpha)}, \cV_1^{(\alpha)}, \dotsc$ is a sequence of covers of $X$.
Consequently, there is an ordinal number $\alpha<\fc$ such that the sequence $\eseq{\cU}$ is equal to $\cU_0^{(\alpha)}, \cU_1^{(\alpha)}, \dotsc$ .
By~\ref{itm:2men} we have 
\[
\sset{u_\beta,w_\beta}{\beta<\alpha}=\sset{u_\beta,w_\beta}{\beta<\alpha}\cap G_\alpha\sub F_\alpha=\Un_{n\in\w}\Un\cF_n^{(\alpha)}
\]
and by~\ref{itm:3men} we have 
\[
\sset{u_\beta, w_\beta}{\alpha\leq\beta<\fc}\sub F_\alpha=\Un_{n\in\w}\Un\cF_n^{(\alpha)}.
\]
Then the family $\Un_{n\in\w}\cF_n^{(\alpha)}$ covers $X$, and thus the set $X$ is Menger.

By~\ref{itm:4men} the set $X$ is totally imperfect and by~\ref{itm:5men} we have $\pi[X]=\Cantor$.
\epf

\bcor
Assume that $\cov(\cM)=\fc$.
There is a Menger totally imperfect set which is not $\sgo$.
\ecor

\section{More consequences}

\subsection{On the problem of Nowik and Tsaban}

The set from Theorem~\ref{thm:main} provides a negative answer to the question from Problem~\ref{prb:pmt} under $\cov(\cN)=\fb=\fc$.
On the other hand it is consistent with ZFC that all Hurewicz totally imperfect sets have cardinalities strictly smaller than $\fc$.
By the results of Miller~\cite{mill} and Zakrzewski~\cite{zak}, this holds in the Cohen model, and by the result of Miller~\cite{mill} also in the Sacks model.
In these cases an affirmative answer to Problem~\ref{prb:pmt} follows from Theorem~\ref{thm:nowik}.

\subsection{Everywhere meager and everywhere null}
\label{sec:EM}
A subset $A$ of a Polish space $X$ is \emph{universally meager}~\cite{zak, zak2,bart, bart2} if for any Borel isomorphism between $X$ and a Polish space $Y$, the image of $A$ under the map is meager.
A set $X\sub\Cantor$ is \emph{everywhere meager}~\cite{krasz} if for each set $a\in\roth$ the set $X_a:=\sset{x\restrict a}{x\in X}$ is meager in $2^a$. 
Kraszewski proved that there is a universally meager set in $\Cantor$ which is not everywhere meager~\cite[Theorem~4.1]{krasz}.
We give a consistently negative answer to a weaker version of the following question.

\bprb[{Kraszewski~\cite[Problem~1]{krasz}}]\label{prb:krasz}
Is every PMT subset of $\Cantor$ everywhere meager?
\eprb

A set $X\sub\Cantor$ is \emph{countably perfectly meager} if for any sequence $\eseq{P}$ of perfect subsets of $\Cantor$, there is an $F_\sigma$-set $F\sub\Cantor$ containing $X$ and the sets $P_n\cap F$ are meager in $P_n$ for all natural numbers $n$.
Equivalently, the set $X$ is countably perfectly meager if for each perfect set $P$ and a countable set $Q$ there is an $F_\sigma$-set $F$ containing $X$ such that for every $q\in Q$, the set $F\cap(q + P)$ is meager in $q+P$. 
Every Hurewicz totally imperfect subset of $\Cantor$ is countably perfectly meager~\cite{zakp}.
By the results of Pol--Zakrzewski~\cite{zakp} and Weiss--Zakrzewski~\cite{zakw} there is consistently a universally meager set which is not countably perfectly meager.

\blem\label{lem:em}
If $X\sub\Cantor$ is a set whose continuous image into $\Cantor$ is not meager, then a homeomorphic copy of $X$ in $\Cantor$ is not everywhere meager.
\elem

\bpf
Let $f\colon X\to \Cantor$ be a continuous map, whose image $f[X]$ is not meager in $\Cantor$.
Then the set 
\[
Y:=\sset{(x,f(x))}{x\in X},
\]
the graph of the function $f$, is a homeomorphic copy of $X$ in $\Cantor\x\Cantor$.
Fix disjoint sets $a,b\in\roth$ with $a\cup b=\w$.
Let $\psi_a\colon \Cantor \to 2^a$ and $\psi_b\colon \Cantor\to 2^b$ be homeomorphisms.
Then the map $\psi:=\psi_a\x\psi_b$ is a homeomorphism from $\Cantor\x\Cantor$ to $2^a\x2^b$.
Let $\phi\colon 2^a\x2^b\to \Cantor$ be a homeomorphism such that $\phi(x,y):=x\cup y$ for all $(x,y)\in 2^a\x2^b$.
Then the set 
\[
Z:=(\phi\circ\psi)[Y]
\]
is a homeomorphic copy of $X$ such that the set $Z_b=\psi_b[f[X]]$ is not meager in $2^b$.
\epf

The following result answers Problem~\ref{prb:krasz} in which PMT sets are replaced by countably perfectly meager sets in the negative.

\bcor
Assume that $\cof(\cN)=\fb=\fc$.
Then there is a totally imperfect Hurewicz (and thus countably perfectly meager) subset of $\Cantor$ which is not everywhere meager.
\ecor

\bpf
Let $X$ be a set from Theorem~\ref{thm:main}.
By Lemma~\ref{lem:em}, there is a homeomorphic copy $Y$ of $X$ in $\Cantor$ which is not everywhere meager.
Since $Y$ is totally imperfect and Hurewicz, by the result of Pol--Zakrzewski~\cite{zakp}, it is countably perfectly meager.
\epf

According to a characterization of universally meager sets due to Zakrzewski~\cite[Theorem~2.1]{zak}, a subset of a Polish space is universally meager if and only if every Borel isomorphic image of this set into a Polish space is meager.
It follows that universally meager sets are preserved by homeomorphisms.
By the result of Pol and Zakrzewski~\cite[Theorem~3.1.3]{zakp} consistently countably perfectly meager sets are not preserved by homeomorphisms.
We do not know whether of everywhere meager sets are preserved by homeomorphisms.

\section{Comments and open problems}

\subsection{Groupability and considered properties}

Considered here properties as Hurewicz and Menger have also equivalent formulations using another types of open covers of spaces.
Let $\cU$ be a countable cover of a set $X$.
The family $\cU$ is an \emph{$\w$-cover} of $X$, if $X\notin U$ and each finite subset of $U$ is contained in a set from $\cU$. 
The family $\cU$ is \emph{groupable}, if there is a partition $\sset{\cF_n}{n\in\w}$ of $\cU$ into pairwise disjoint finite sets such that the family $\sset{\Un\cF_n}{n\in\w}$ is a $\gamma$-cover of $X$.
A set $X$ satsifies $\sfin(\Omega,\Op)$ if for any sequence $\eseq{\cU}$ of open $\w$-covers of $X$, there are finite sets $\cF_0\sub\cU_0,\cF_1\sub\cU_1,\dotsc$ such that the family $\Un_{n\in\w}\cF_n$ is a cover of $X$.
A set $X$ satsifies $\sfin(\Omega,\Op^\mathrm{gp})$ if for any sequence $\eseq{\cU}$ of open $\w$-covers of $X$, there are finite sets $\cF_0\sub\cU_0,\cF_1\sub\cU_1,\dotsc$ such that the family $\Un_{n\in\w}\cF_n$ is a groupable cover of $X$.
By a result of Scheepers~\cite{coc1} the Menger property is equivalent to $\sfin(\Omega,\Op)$ and by the result of Kocinac--Scheepers~\cite[Theorem~14]{grp} the Hurewicz property is equivalent to $\sfin(\Omega,\Op^\mathrm{gp})$.
Our paper provides another than already known examples of $\sfin(\Omega,\Op)$ totally imperfect sets which are not $\sfin(\Omega,\Op^\mathrm{gp})$. 
Using dichotomic proof, Chaber--Pol~\cite{chp} proved that there is such a set in ZFC.
A unfiorm ZFC construction of such a set was provied by Tsaban--Zdomskyy~\cite{sfh}.

\subsection{The property $\un(\Gamma,\Gamma)$}
A set $X$ satisfies $\un(\Gamma,\Gamma)$, if for any sequence $\eseq{\cU}$ of open $\gamma$-covers of $X$, there are sets $\cF_0\sub\cU_0,\cF_1\sub\cU_1,\dotsc$ such that $\card{\cF_n}=n$ for all $k$'s and the family $\sset{\Un\cF_n}{n\in\w}$ is a $\gamma$-cover of $X$.
Let $k$ be a positive natural number.
A set $X$ satisfies $\unn(\Gamma,\Gamma)$, if for any sequence $\eseq{\cU}$ of open $\gamma$-covers of $X$, there are sets $\cF_0\sub\cU_0,\cF_1\sub\cU_1,\dotsc$ such that $\card{\cF_n}=k$ for all $n$'s and the family $\sset{\Un\cF_n}{n\in\w}$ is a $\gamma$-cover of $X$.
This property was investigated by Tsaban~\cite{thebook} in the context of sets with the Hurewicz property and the property $\sgg$. 
A set $X$ satisfies $\sgg$ if for any sequence $\eseq{\cU}$ of open $\gamma$-covers of $X$, there are sets $U_0\in\cU_0, U_1\in\cU_1,\dotsc$ such that the family $\sset{U_n}{n\in\w}$ is a $\gamma$-cover of $X$.
We have the following implications

\begin{figure}[H]
\begin{center}
\begin{tikzcd}[ampersand replacement=\&,column sep=1cm]
{}\& {}\& \text{Hurewicz}\&{}\&{}\\
\sgg\arrow{r}\&  \unn(\Gamma,\Gamma)\arrow{r}\&\un(\Gamma,\Gamma)\arrow{u}\arrow{r}\& \sgo\arrow{r}\&\text{totally imperfect}.
\end{tikzcd}
\end{center}
\end{figure}

Tsaban asked the following question.

\bprb[{\cite[Problem~3.9]{thebook}}]
Is there a totally imperfect Hurewicz set which is not $\un(\Gamma,\Gamma)$?
\eprb

Our example constructed in Theorem~\ref{thm:main} provides a negative answer to the above question under $\cov(\cN)=\fb=\fc$.
Assume that $\fb=fc$.
In a very recent paper~\cite{omi}, 
using different techniques Tsaban proved that there is a Hurewicz set satisfying $\sgo$ (and thus, it is totally imperfect) which is not $\un(\Gamma,\Gamma)$~~\cite[Theorem~10.6]{omi}.
Tsbaban also proved that there is a set satisfying $\un(\Gamma,\Gamma)$ but not $\unn(\Gamma,\Gamma)$ for any positive natural number $k$~\cite[Theorem~9.9]{omi} and also that for  each $k$ there is a set satisfying $\mathsf{U}_{k+1}(\Gamma,\Gamma)$ but not $\unn(\Gamma,\Gamma)$\cite[Theorem~10.3]{omi}
There are also another results of Liu--He--Zhang in a similar spirit, if we assume that the Continuum Hypothesis holds~\cite{lhz}.

\subsection{Separation of properties}

\bprb
Is there in ZFC a totally imperfect Menger set which is not perfectly meager? Is there in ZFC a perfectly meager set which is not Menger\footnote{In a recent work of V. Haberl with the first and third named authors it is shown that every perfectly meager set has size $\omega_1$ in the Miller model, and hence all such sets are Menger there.}?
\eprb

\bprb
Is there a Menger set whose all continuous images into $\Cantor$ are totally imperfect and which does not satisfy $\sgo$?
\eprb

\subsection{Around PMT sets}

\bprb
 Does the existence of a totally imperfect Hurewicz set of size $\fc$ imply that there is a totally imperfect Hurewicz set which is not PMT?
\eprb

\bprb
 Let $X\sub\Cantor$.
 Are the following statements equivalent?
 \be
 \item The set $X$ is Hurewicz and totally imperfect which can be mapped continuously onto $\Cantor$.
 \item The set $X$ is not PMT?
 \ee
\eprb

\section*{Acknowledgments}
We would like to thank to the anonymous referee for his/her suggestions and comments which help us to increase substantially the quality of this paper.
We also would like to thank to Professor Piotr Zakrzewski and Professor Roman Pol for fruitful discussions about obtained in the paper results.

\end{document}